# Probability Bracket Notation, Wick-Matsubara Relation, Density Operators, and Microscopic Probability Modeling

Xing M. Wang [1]

## Abstract

Our previous paper introduced the probability bracket notation (PBN) following the Dirac vector bracket notation (VBN). We demonstrated that, under a special Wick rotation (imaginary time), the stationary Schrödinger equation in Hilbert space transforms into the master equation of a microscopic probabilistic process (MPP) in the probability space. In this article, we first analyze the MPP of a single-particle system, showing that the energy expectation of the MPP eventually converges to the lowest energy level based on its initial conditions. At the same time, its von Neumann entropy ultimately approaches zero. Next, we explore the MPP for a quantum system of identical particles in the Fock space. By connecting time with temperature through the Wick-Matsubara relation, we recover the expected occupation number of particles and the grand partition function in quantum statistics. Furthermore, we reproduce the internal energy of an ideal gas in thermodynamics using this relation. To address the entropy issue and adapt the PBN to statistical research topics in the literature, we express the density operators and von Neumann entropy in probability space using the PBN framework. Integrating the Wick-Matsubara relation with the PBN may provide a new perspective on microscopic probability modeling, thermodynamics, and quantum statistics.

[1] Sherman Visual Lab, Sunnyvale, CA, USA, xmwang@shermanlab.com

## 1. Introduction

The Dirac vector probability notation (VBN, see [1], p.96) is a powerful tool to manipulate vectors in the Hilbert space. The main advantage of the VBN is that many formulas in QM can be presented in a symbolic and abstract form, independent of the state expansion or the basis selection, which, when needed, can be easily done by inserting an identity operator. Inspired by the success of VBN, we proposed the Probability Bracket Notation (PBN) [2]. Below is a brief introduction to both notations.

The V-bracket in the VBN is an inner product of two vectors in the Hilbert space.

$$\text{V-bracket: } \langle\psi_A|\psi_B\rangle \quad\Rightarrow\quad \text{V-bra}: \ \langle\psi_A|, \quad \text{V-ket}: \ |\psi_B\rangle \tag{1.1}$$

$$\text{where}: \ \langle\psi_A|\psi_B\rangle = \langle\psi_B|\psi_A\rangle^*, \quad \langle\psi_A| = |\psi_A\rangle^\dagger \tag{1.2}$$

V-basis and V-identity:

$$\text{A discrete basis:} \quad \hat{H}|\varepsilon_i\rangle = \varepsilon_i|\varepsilon_i\rangle, \quad \langle\varepsilon_i|\varepsilon_j\rangle = \delta_{ij}, \quad \hat{I}_H = \sum_i |\varepsilon_i\rangle\langle\varepsilon_i| \tag{1.3}$$

$$\text{A continuous basis:} \quad \hat{x}|x\rangle = x|x\rangle, \quad \langle x|x'\rangle = \delta_{x,x'}, \quad \hat{I}_X = \int dx|x\rangle\langle x| \tag{1.4}$$

The normalization of a system V-ket $|\Psi\rangle$ using a continuous V-basis:

$$1 = \langle\Psi|\Psi\rangle = \langle\Psi|\hat{I}_X|\Psi\rangle = \int dx\langle\Psi|x\rangle\langle x|\Psi\rangle = \int dx|\Psi(x)|^2 \tag{1.5}$$

The expectation of an operator $H$ with a discrete V-basis:

$$\langle H\rangle \equiv \bar{H} \equiv \langle\Psi|\hat{H}|\Psi\rangle = \sum_i \langle\Psi|\hat{H}|\varepsilon_i\rangle\langle\varepsilon_i|\Psi\rangle = \sum_i \varepsilon_i |c_i|^2 \tag{1.6}$$

The P-bracket in the PBN [2] is a conditional probability in the probability space.

$$\text{P-bracket: } P(\varphi_A|\varphi_B) \quad\Rightarrow\quad \text{P-bra}: \ P(\varphi_A|, \quad \text{P-ket}: \ |\varphi_B) \tag{1.7}$$

$$\text{where}: \ P(\varphi_A|\varphi_B) \neq P(\varphi_B|\varphi_A)^*, \quad P(\varphi_A| \neq |\varphi_A)^\dagger, \quad P(A|B) = 1, \text{if } B \subseteq A \tag{1.8}$$

P-basis and P-identity:

$$\text{A discrete case:} \quad \hat{H}|\varepsilon_i) = \varepsilon_i|\varepsilon_i), \quad P(\varepsilon_i|\varepsilon_j) = \delta_{ij}, \quad \hat{I}_H = \sum_i |\varepsilon_i)P(\varepsilon_i| \tag{1.9}$$

$$\text{A continuous case:} \quad X|x) = x|x), \quad P(x|x') = \delta_{x,x'}, \quad \hat{I}_X = \int_{x\in\Omega} dx|x)P(x| \tag{1.10}$$

The normalization of a system P-ket with a continuous P-basis:

$$1 = P(\Omega | \Omega) = P(\Omega | I_X | \Omega) = \int_{x\in\Omega} dxP(\Omega | x)P(x | \Omega) = \int_{x\in\Omega} dxP(x | \Omega) = \int_{x\in\Omega} dxP(x) \quad (1.11)$$

The expectation of an operator *H* with a discrete P-basis:

$$\langle H \rangle \equiv \bar{H} \equiv P(\Omega | \hat{H} | \Omega) = \sum_i P(\Omega | \hat{H} | \varepsilon_i)P(\varepsilon_i | \Omega) = \sum_i \varepsilon_i P(\varepsilon_i | \Omega) = \sum_i \varepsilon_i P(\varepsilon_i) \quad (1.12)$$

According to the Born statistic interpretation ([1], p.53) of the wave functions, there exists a relation between the system V-ket and the system P-ket:

$$| \Psi_t \rangle = \sum_i c_i(t) | \psi_i \rangle \underset{\text{one to may}}{\overset{\text{many to one}}{\rightleftarrows}} | \Omega_t ) = \sum_i P(i,t) | i); \quad P(i,t) = | c_i(t) |^2 \quad (1.13)$$

Then, for example, the expectation value can be written in both ways:

$$\langle H \rangle_t = \langle \Psi_t | \hat{H} | \Psi_t \rangle = \sum_i \varepsilon_i | c_i(t) |^2 = P(\Omega | \hat{H} | \Omega_t \rangle = \sum_i \varepsilon_i P(i,t) \quad (1.14)$$

Therefore, the PBN is a very convenient tool for handling the microscopic probability models based on Eq. (1.13-14).

In Ref. [2], we discussed several correlations revealed by the special Wick rotation (SWR, $it \rightarrow t$). In particular: the stationary Schrodinger equations in the Hilbert space transform into the master equations of microscopic probabilistic processes (MPPs) in the probability space under the SWR.

In the following section, we will investigate the time evolution of the MPP for a single particle with either a discrete or a continuous energy spectrum. Our findings suggest that an MPP behaves similarly to the freezing process. We will then explore the time evolution of the MPP for the system of many identical particles in Fock space. By relating time to temperature through the Wick-Matsubara relation $it \rightarrow t \rightarrow \hbar / kT$, we can reproduce the expected occupation number of particles and the grand partition function in quantum statistics. Using the relation, we can also recover the internal energy of an ideal gas in thermodynamics. To address the entropy issue in an MPP and adapt the PBN to research subjects in the statistical literature, we will express the density operators and von Neumann entropy in probability space using the PBN.

## 2. The Microscopic Probabilistic Processes of a Single Particle

By studying the path integrals with PBN [2], we demonstrated that the stationary Schrödinger equation in Hilbert space naturally transforms into the master equation in probability space under a special Wick rotation (SWR). The master equation describes the time evolution of an induced microscopic diffusion, which we call a macroscopic probabilistic process (MPP).

The **Special Wick Rotation** (SWR) is defined by [2]:

$$\text{SWR: } i\,t \to t, \quad |\psi(t)\rangle \to |\Omega_t), \quad \langle x_b, t_b \,|\, x_a, t_a\rangle \to P(x_b, t_b \,|\, x_a, t_a) \tag{2.1}$$

Under the SWR, a stationary Schrodinger equation becomes a master equation of an MPP:

$$\frac{\partial}{i\,\partial t}|\Psi(t)\rangle = \frac{-1}{\hbar}\hat{H}\,|\psi(t)\rangle \overset{it\to t}{\to} \frac{\partial}{\partial t}|\Omega_t) = \frac{-1}{\hbar}\hat{H}\,|\Omega_t), \quad \partial_t \hat{H} = 0 \tag{2.2}$$

In Ref. [2], we discussed the MPP of a single particle and its possible application to certain systems with non-Hermitian Hamiltonians [3]. In this section, we will discuss the behavior of the MPPs in detail.

**2.1. A Particle with a Discrete Energy Spectrum**: For the system of one particle with a discrete energy spectrum, the time evolution of the energy eigenstates exhibits the following correlation under the SWR (Eq. (115), [2]):

$$|\psi_k(t)\rangle = e^{-i\varepsilon_k t/\hbar}\,|\psi_k\rangle \overset{it\to t}{\to} |\psi_k(t)) = e^{-\varepsilon_k t/\hbar}\,|\psi_k) \tag{2.3}$$

Suppose that the MPP has a normalized initial system P-ket:

$$|\Omega_0) = \sum\nolimits_{k\in B} \eta_k\,|\psi_k), \quad P(\Omega\,|\,\Omega_0) = \sum\nolimits_{k\in B} \eta_k = 1 \tag{2.4}$$

It leads to the following unnormalized time-dependent system P-ket,

$$|\Omega_t) = \sum\nolimits_{k\in B} \eta_k e^{-\varepsilon_k t/\hbar}\,|\psi_k), \quad P(\Omega\,|\,\Omega_t) = \sum\nolimits_{k\in B} \eta_k e^{-\varepsilon_k t/\hbar} \neq 1 \quad \text{for } t > 0 \tag{2.5}$$

Because it is unnormalized, the expectation of energy has the expression:

$$\langle \hat{H}\rangle_t = P(\Omega\,|\,\hat{H}\,|\,\Omega_t)\,/\,P(\Omega\,|\,\Omega_t) \tag{2.6}$$

Due to the exponential time dependence in Eq. (2.5), the basis P-ket of the lowest energy level ultimately dominates the system P-ket. Let $\varepsilon_s, s \in B$ be the lowest energy, we have

$$|\Omega_t) = \sum\nolimits_{k\in B} \eta_k e^{-\varepsilon_k t/\hbar}\,|\psi_k) \overset{t\to\infty}{\to} \eta_s e^{-\varepsilon_s t/\hbar}\,|\psi_s) \tag{2.7}$$

Therefore, the MPP behaves like a process that will eventually freeze to its lowest energy level, defined in its initial condition (2.4) (see the discussion after Eq. (117), [2]):

$$\langle \hat{H}\rangle_t = P(\Omega\,|\,\hat{H}\,|\,\Omega_t)\,/\,P(\Omega\,|\,\Omega_t) = \sum\nolimits_{k\in B} \eta_k\,\mathrm{e}^{-\omega_k t}\,\varepsilon_k \,/\, \sum\nolimits_{k\in B} \eta_k\,\mathrm{e}^{-\omega_k t} \overset{t\to\infty}{\to} \varepsilon_s \tag{2.8}$$

The system reaches a pure state with zero von Neumann entropy [4] and the MPP resembles a freezing process.

**2.2. A particle with a Continuous Energy Spectrum**: Let us now consider a free particle that has a uniform probability distribution in momentum $p$. According to Eq. (2.3), its system P-ket is represented as:

$$|\Omega_t) = \int_p dp\, e^{-\varepsilon_p t/\hbar}\, | p) = \int_p dp\, e^{-\varepsilon_p \beta}\, | p), \quad \varepsilon_p \equiv p^2/2m, \beta \equiv t/\hbar \tag{2.9}$$

Because it is unnormalized, according to Eq. (2.6), the expected energy reads:

$$\langle \hat{H} \rangle = P(\Omega | \hat{H} | \Omega_t) / P(\Omega | \Omega_t) = \int_p dp \varepsilon_p e^{-\varepsilon_p \beta} / \int_p dp e^{-\varepsilon_p \beta} = -\partial_\beta \ln z_\beta \tag{2.10}$$

$$z_\beta \equiv \int_p dp\, e^{-\varepsilon_p \beta} \tag{2.11}$$

It can be easily calculated (see [5], §10.4):

$$z_\beta = \int dp\, e^{-\beta\varepsilon} = \int dp\, e^{-\beta p^2/2m} = \sqrt{\pi\beta/2m} \tag{2.12}$$

$$\langle \hat{H} \rangle = -\partial_\beta \ln z_\beta = \frac{1}{2\beta} = \frac{\hbar}{2t} \tag{2.13}$$

The explicit time dependency indicates that, at $t = 0$, the system has infinite expected energy due to its uniformly distributed momentum. As time progresses towards infinity, the system's energy decreases and approaches zero. This scenario resembles a microscopic freezing process, where the energy level eventually falls to the lowest ($\varepsilon_p = 0$ in the continuous $p$-spectrum).

It is reasonable to question whether the time in Eq. (2.5) or (2.13), initially conceived as imaginary time, might be connected to temperature. The answer seems to be YES.

## 3. Many-Particle Systems and the Wick-Matsubara Relation

3.1. **The MPP of a Many-Particle System**: For a system of many identical and weakly interacting particles, we utilize the following P-basis of occupation numbers, mapped from the V-basis in the Fock space (Eq. (28), [2]):

$$\hat{n}_i | n_1, n_2, \ldots) = \hat{n}_i \prod_k | n_k) \equiv n_i | \vec{n}) \tag{3.1}$$

$$P(\vec{n}' | \vec{n}) = \prod_i P(n_i' | n_i) \equiv \delta_{\vec{n}',\vec{n}}, \quad \sum_{\vec{n}} | \vec{n}) P(\vec{n} | = \prod_i \sum_{n_i} | n_i) P(n_i | = \hat{I}_{\vec{n}} \tag{3.2}$$

The system has the following total Hamiltonian $\boldsymbol{H}$ (Eq. (1.3), [6]):

$$\hat{\boldsymbol{H}} = \hat{H} - \mu\hat{N} = \sum_i \sum_{n_i} \hat{n}_i(\varepsilon_i - \mu) \tag{3.3}$$

Here $\hat{H}$, $\varepsilon_i$ and $\mu$ are the conventional dynamic Hamiltonian, the dynamic energy per particle at the $i^{\text{th}}$ level, and the chemical potential per particle. Based on Eq. (2.3) and (3.3), the time evolution of the basis P-ket at the $i^{\text{th}}$ level is:

$$\partial_t |n_i, t) = -[\hat{n}_i(\varepsilon_i - \mu)/\hbar]\,|n_i, t), \quad |n_i, t) = e^{-n_i(\varepsilon_i - \mu)t/\hbar}\,|n_i) \tag{3.4}$$

We assume that at $t = 0$, the probability of the system is uniformly distributed:

$$|\Omega_{i,0}) \equiv |\Omega_i) = \sum_{n_i} |n_i), \quad |\Omega_{i,t}) = \sum_{n_i} e^{-n_i(\varepsilon_i - \mu)t/\hbar}\,|n_i) \tag{3.5}$$

It is unnormalized and the expected occupation number of particles at the $i^{\text{th}}$ level reads:

$$n(\varepsilon_i) \equiv \langle \hat{n}_i \rangle = \frac{\sum_{n_i} n_i e^{-(\varepsilon_i - \mu)n_i t/\hbar}}{\sum_{n_i} e^{-(\varepsilon_i - \mu)n_i t/\hbar}} \tag{3.6}$$

The expected total number of particles can be obtained:

$$\langle \hat{N} \rangle = \sum_i \langle \hat{n}_i \rangle == \sum_i \left\{ \sum_{n_i} n_i e^{-(\varepsilon_i - \mu)n_i t/\hbar} \Big/ \sum_{n_i} e^{-(\varepsilon_i - \mu)n_i t/\hbar} \right\} \tag{3.7}$$

It exhibits the same freezing behavior as the single particle with a discrete energy spectrum, as indicated in Eq. (2.8). However, Eq. (3.7) bears a close resemblance to the distribution function found in quantum statistics (11.2-6, [5]):

$$\langle \hat{N} \rangle = \sum_i \left\{ \sum_{n_i} n_i e^{-(\varepsilon_i - \mu)n_i/kT} \Big/ \sum_{n_i} e^{-n_i(\varepsilon_i - \mu)/kT} \right\} \equiv \sum_i \sum_{n_i} n_i e^{-(\varepsilon_i - \mu)n_i/kT} / Z_i \tag{3.8}$$

Here $Z_i$ is the partition function of the $i^{\text{th}}$ energy level.

**3.2. The Wick-Matsubara Relation**: By comparing Eq. (3.7) with (3.8), we observe that the imaginary time in the Wick rotation corresponds to absolute temperature in the Matsubara formalism (see §1.4, [6]). We designate this as:

**The Wick-Matsubara Relation**: $$it \xrightarrow{\text{Wick}} t \xrightarrow{\text{Matsubara}} \hbar/kT \tag{3.9}$$

For fermions, $n_i \in \{0, 1\}$ (the Fermi-Dirac distribution), all lowest energy states ($\varepsilon \leq \varepsilon_F$) are occupied when $T \to 0$ (or $t \to \infty$). For bosons, $n_i \geq 0$ (the Bose-Einstein distribution), when $T \to 0$ (or $t \to \infty$), all boson go to the ground energy level (the Bose condensation).

At zero temperature, the systems of either bosons or fermions reach a frozen state, characterized by the lowest total energy and zero entropy (no chaos).

Using Eq. (3.9), Eq. (2.13) becomes:

$$\langle \hat{H} \rangle_t = \frac{\hbar}{2t} \rightarrow \langle \hat{H} \rangle = \frac{kT}{2} \tag{3.10}$$

Eq. (3.10) represents the expected energy of one particle in one dimension of an ideal monoatomic gas ($\langle \varepsilon \rangle = kT/2$). For a system of $N$ non-interacting identical particles in 3D, we recover its internal energy in thermodynamics (§10.7, [5]):

$$U = 3N\langle \hat{H} \rangle = \frac{3NkT}{2} \tag{3.11}$$

The grand partition function can be expressed in the VBN as (§11.2-6, [5]; page 37, [6]):

$$Z_G \equiv \sum\nolimits_{\vec{n}} \langle \vec{n} \,|\, e^{-(\hat{H}-\mu\hat{N})/kT} \,|\, \vec{n} \rangle = \prod\nolimits_{i=1}^{\infty} \sum\nolimits_{n_i} e^{-(\varepsilon_i-\mu)n_i/kT} = \prod\nolimits_{i=1}^{\infty} Z_i \tag{3.12}$$

We can write it concisely in the PBN and it is easy to expand:

$$Z_G = P(\Omega \,|\, e^{-(\hat{H}-\mu\hat{N})/kT} \,|\, \Omega) = P(\Omega \,|\, e^{-(\hat{H}-\mu\hat{N})/kT} \hat{I}_{\vec{n}} \,|\, \Omega) = \prod\nolimits_{i=1}^{\infty} \sum\nolimits_{n_i} e^{-(\varepsilon_i-\mu)\hat{n}_i/kT} \tag{3.13}$$

Eq. (3.6) can also be rewritten in a compact form:

$$\langle \hat{n}_i \rangle = \frac{P(\Omega \,|\, \hat{n}_i \,|\, \Omega_t)}{P(\Omega \,|\, \Omega_t)} = \frac{P(\Omega \,|\, \hat{n}_i e^{-\hat{\boldsymbol{H}} t/\hbar} \,|\, \Omega)}{P(\Omega \,|\, e^{-\hat{\boldsymbol{H}} t/\hbar} \,|\, \Omega)} = \frac{P(\Omega \,|\, \hat{n}_i e^{-\hat{\boldsymbol{H}} t/\hbar} \hat{I}_{\vec{n}} \,|\, \Omega)}{P(\Omega \,|\, e^{-\hat{\boldsymbol{H}} t/\hbar} \hat{I}_{\vec{n}} \,|\, \Omega)} = \frac{\sum\nolimits_{n_i} n_i e^{-(\varepsilon_i-\mu)n_i t/\hbar}}{\sum\nolimits_{n_i} e^{-(\varepsilon_i-\mu)n_i t/\hbar}} \tag{3.14}$$

The freezing behavior of an MPP implies that the system appears to be in contact with an enormous reservoir at an absolute zero temperature. In reality, the system will maintain the temperature $T$ of its environment (or the reservoir in contact). Consequently, its time clock (imaginary time, not real!) will quickly halt at $t = \hbar/kT$. For example, if the reservoir is at room temperature, the freezing process will instantly cease at the time:

$$t \sim \frac{\hbar}{kT} \sim \frac{1.05\times10^{-27}}{1.38\times10^{-16}\times300} \sim 2.54\times10^{-14}\,(\text{s}) \tag{3.15}$$

Since the master equation of an MPP describes a probabilistic process, even though Eq. (2.5) and Eq. (2.13) pertain to a single particle, the Wick-Matsubara relation can be applied, associating their time with temperature.

## 4. Density Operators and the von Neumann Entropy in the PBN

Up to this point, we have not utilized the density operators in the probability space because the system P-ket already contains all the information for expectations, as shown in Eq. (2.6). However, we do need these operators to explore the von Neumann entropy in MPPs and to connect the PBN with relevant statistical research topics in the literature [6].

**4.1. Density Operators in Hilbert Space**: The density operator [7] (or matrix) $\rho$ is defined in Hilbert space $\mathcal{H}^N$, having the following properties:

$$\mathrm{Tr}\rho = 1, \quad \rho = \rho^\dagger, \quad \langle \varphi | \rho | \varphi \rangle \geq 0 \ \ \forall \ |\varphi\rangle \in \mathcal{H}^N \tag{4.1}$$

A density matrix $\rho$ describes a pure state [8] if and only if $\rho^2 = \rho$, i.e. the state is idempotent; it describes a mixed state [9] if $\rho^2 \neq \rho$. The purity of a state is defined as $\mathrm{Tr}\rho^2$. For a pure state, $\mathrm{Tr}\rho^2 = 1$, and for a mixed state, $\mathrm{Tr}\rho^2 < 1$.

Let $\eta_k$, $k \in \{1, 2, ..., \}$, be the eigenvalues of the density matrix, then we have the following expressions:

$$\rho | \eta_k \rangle = \eta_k \, | \eta_k \rangle, \quad \langle \eta_i | \eta_k \rangle = \delta_{i,k}, \quad \rho = \sum\nolimits_k | \eta_k \rangle \, \eta_k \langle \eta_k |, \quad \mathrm{Tr}\rho = \sum\nolimits_k \eta_k = 1 \tag{4.2}$$

**4.2. Density Operators in the Probability Space for MPPs**: Using the PBN, we can define the density operator with the system P-ket and its basis, similar to Eq. (4.2) in Hilbert space. We define the time-dependent density operator for an MPP as follows:

$$\rho_t \equiv \sum\nolimits_k [\eta_k e^{-\varepsilon_k t/\hbar} \, | \psi_k ) \, P(\psi_k \, |] / \sum\nolimits_k \eta_k e^{-\varepsilon_k t/\hbar}, \quad \langle \hat{H} \rangle_t = \mathrm{Tr}[\hat{H} \rho_t] \tag{4.3}$$

Here, we assume that the Hamiltonian of the MPP has a discrete spectrum as in Eq. (2.5):

$$| \Omega_t ) = \sum\nolimits_{k \in B} \eta_k e^{-\varepsilon_k t/\hbar} \, | \psi_k ), \quad P(\Omega \, | \, \Omega_t) = \sum\nolimits_{k \in B} \eta_k e^{-\varepsilon_k t/\hbar} \neq 1 \quad \text{for } t > 0 \tag{4.4}$$

The time-dependent density operator can also be written as:

$$\rho_t = \sum\nolimits_k [\eta_k(t) \, | \, \psi_k ) P(\psi_k \, |] / \sum\nolimits_k \eta_k(t), \quad \eta_k(t) = \eta_k e^{-\varepsilon_k t/\hbar} \tag{4.5}$$

The expectation of the Hamiltonian $H$ is given by:

$$\langle \hat{H} \rangle_t = \mathrm{Tr}[\hat{H} \rho_t] = \sum\nolimits_k \varepsilon_k \eta_k(t) / \sum\nolimits_k \eta_k(t) \tag{4.6}$$

Similarly, for the system of a free particle in Eq. (2.9), the density operator is:

$$\rho_t = \frac{\int dp \, e^{-\varepsilon_p t/\hbar} \, | \, p) P(p \, |}{\int dp \, e^{-\varepsilon_p t/\hbar}} = \frac{e^{-\hat{H} t/\hbar} \hat{I}_p}{\mathrm{Tr}[e^{-\hat{H} t/\hbar} \hat{I}_p]} = \frac{e^{-\beta \hat{H}}}{z_\beta}, \quad \varepsilon_p = p^2 / 2m, \ \ \beta = t / \hbar \tag{4.7}$$

$$\langle \hat{H} \rangle_t = \mathrm{Tr}[\hat{H}\rho_t] == \frac{\mathrm{Tr}[\hat{H}e^{-\beta\hat{H}}]}{\mathrm{Tr}[e^{-\beta\hat{H}}]} = \frac{\int dp\, \varepsilon_p e^{-\varepsilon_p \beta}}{z_\beta} = -\partial_\beta \ln z_\beta \tag{4.8}$$

For the system of many particles in Section 2.2, the density operator reads:

$$\rho_\beta = \frac{\sum_{\vec{n}} e^{-n_i(\varepsilon_i-\mu)t/\hbar} |\vec{n})P(\hat{n}|}{\sum_{\vec{n}} e^{-n_i(\varepsilon_i-\mu)t/\hbar}} = \frac{e^{-\hat{\boldsymbol{H}}\beta}\hat{I}_{\vec{n}}}{\mathrm{Tr}[e^{-\hat{\boldsymbol{H}}\beta}\hat{I}_{\vec{n}}]} = \frac{e^{-\hat{\boldsymbol{H}}\beta}}{Z_G(\beta)}, \quad \beta \equiv \frac{t}{\hbar} \tag{4.9}$$

The expected occupation number of particles at the $i^{\mathrm{th}}$ level in Eq. (3.14) becomes:

$$\langle \hat{n}_i \rangle_t = \mathrm{Tr}(\hat{n}_i \rho_t) = \frac{\mathrm{Tr}(\hat{n}_i e^{-\boldsymbol{H}t/\hbar})}{Z_G(\beta)} = \frac{\sum_{n_i} n_i\, e^{-n_i(\varepsilon_i-\mu)t/\hbar}}{\sum_{n_i} e^{-n_i(\varepsilon_i-\mu)t/\hbar}} \tag{4.10}$$

The grand partition function in Eq. (3.13) can be written as (here $\beta = 1/kT$):

$$Z_G = \mathrm{Tr}\rho_\beta = \mathrm{Tr}[e^{-(\hat{H}-\mu\hat{N})\beta}\hat{I}_{\vec{n}}] = \prod_{i=1}^{\infty} \sum_{n_i} e^{-(\varepsilon_i-\mu)\hat{n}_i/kT} = \prod_{i=1}^{\infty} Z_i \tag{4.11}$$

**4.3. Quantum Entropy** (QE): For a quantum system, the quantum entropy (QE) or the von Neumann entropy [10] is defined by:

$$S(\rho) = -\mathrm{Tr}\,(\rho \ln \rho) \tag{4.12}$$

It can be shown that:

$$S(\rho) \ge 0; \quad S(\rho) = 0 \text{ if only if } \rho \text{ describes a pure state.} \tag{4.13}$$

The maximally mixed state in $N$-dimensional Hilbert space $\mathcal{H}^N$ has a density operator proportional to the $N$-dimensional identity and has the maximal entropy [10]:

$$\rho_{\max} = \frac{1}{N}\hat{I}_N, \quad S(\rho_{\max}) = \ln N \tag{4.14}$$

It corresponds to the uniform distribution. Therefore, the distributions in Eq. (2.9) and (3.5) produce the maximal entropy at the initial time.

In the orthonormal basis of $H$ in Eq. (4.5), the QE is given by:

$$S(\rho_t) = -\mathrm{Tr}\,(\rho_t \ln \rho_t) = -\sum_k \overline{\eta}_k(t) \ln \overline{\eta}_k(t), \quad \overline{\eta}_k(t) \equiv \eta_k(t) / \sum_k \eta_k(t) \tag{4.15}$$

As described in Eq. (2.7), when time approaches infinity, we observe that:

$$|\Omega_t) = \sum_{k\in B} \eta_k e^{-\varepsilon_k t/\hbar} |\psi_k) \overset{t\to\infty}{\to} \eta_s e^{-\varepsilon_s t/\hbar} |\psi_s) \tag{4.16}$$

$$\rho_t = \frac{\sum_{k\in B} \eta_k(t) |\psi_k) P(\psi_k|}{\sum_{k\in B} \eta_k(t)} \overset{t\to\infty}{\to} \frac{\eta_s(t) |\psi_s) P(\psi_s|}{\eta_s(t)} = |\psi_s) P(\psi_s| \tag{4.17}$$

The system is in a pure state, and its QE approaches zero. That was why we stated that an MPP system becomes frozen at its lowest energy and with zero entropy if time approaches infinity. Now, we understand that because of the Wick-Matsubara relation, the (imaginary) time of an MPP will be synchronized with its environmental temperature.

Our expressions are comparable with those in the literature. For example, if we set $t/\hbar = 1/kT = \beta$ and $\eta_k = 1$ (uniform distribution) in our expressions of $\rho_t$ in Eq. (4.3), we obtained equations similar to Eq. (1.1) and (1.11) in Ref [6]:

$$\rho_\beta = \frac{\sum_k e^{-\beta \hat{H}} |\varepsilon_k)\, P(\varepsilon_k|}{\sum_k e^{-\beta\varepsilon_k}} = \frac{e^{-\beta\hat{H}} \hat{I}_\varepsilon}{Z_\beta} = \frac{e^{-\beta\hat{H}}}{Z_\beta} \tag{4.18}$$

$$\langle \hat{A} \rangle = \mathrm{Tr}[\hat{A}\, \rho_\beta] = \frac{\mathrm{Tr}[\hat{A}\, e^{-\beta\hat{H}}]}{Z_\beta} \tag{4.19}$$

## 5. Summary and Discussion

After reviewing the core concepts of the probability bracket notation (PBN), we investigated the implications of one of the correlations revealed by the special Wick rotation (SWR): specifically, under the SWR, a stationary Schrödinger Equation in the Hilbert space naturally transforms into the master equation of a microscopic probabilistic process (MPP) in the probability space.

We demonstrated that the expected energy of the MPPs of a single particle approaches its lowest level, and its entropy eventually vanishes as time goes to infinity. This indicates that the MPP system is an open system, interacting with a vast reservoir at an absolute zero temperature. Next, we investigated the MPP within a system of identical particles in Fock space. We were able to recover the fundamental formulas in thermodynamics and quantum statistics by applying the Wick-Matsubara relation $it \to t \to \hbar / kT$. Based on this relation, the MPP system will take the environmental temperature *T,* and its (imaginary) time clock will quickly synchronize with that temperature ($t = \hbar/kT$). Finally, to address the entropy issue and align our presentations with the statistical research topics in the literature, we express the density operators and the von Neumann entropy in the probability space using the PBN.

When combined with the PBN, the Wick-Matsubara relation may offer a novel approach to microscopic probability modeling, thermodynamics, and quantum statistics. The process likely occurs in three steps: first, make an SWR ($it\rightarrow t$), transforming a stationary Schrodinger equation into the master equation of an MPP using the PBN; next, solve the master equation with an appropriate initial condition; finally, substitute time $t$ with $\hbar/kT$.

Certainly, further investigations are necessary to validate the PBN's consistency, usefulness, and limitations.

## Abbreviations

| | |
|---|---|
| MPP | Microscopic Probabilistic Process |
| PBN | Probability Bracket Notation |
| P-basis | Probability-basis |
| P-bra | Probability-bra |
| P-identity | Probability-identity |
| P-ket | Probability-ket |
| QE | Quantum Entropy |
| SWR | Special Wick Rotation |
| VBN | Vector Bracket Notation |
| V-basis | Vector-basis |